\documentclass[11pt]{amsart}
\usepackage[margin=1in]{geometry}
\usepackage{amsmath,amssymb,mathtools}
\usepackage[T1]{fontenc}
\usepackage{lmodern,microtype}
\usepackage[colorlinks=true,linkcolor=blue,citecolor=blue,urlcolor=blue]{hyperref}
\hypersetup{pdftitle={Noncanonical congruences for the smallest parts function},
 pdfauthor={Larry Rolen},
 pdfsubject={Research draft: extension of the spt noncongruence theorem}}
\newtheorem{theorem}{Theorem}[section]
\newtheorem{corollary}[theorem]{Corollary}
\newtheorem{lemma}[theorem]{Lemma}
\newtheorem{proposition}[theorem]{Proposition}
\numberwithin{equation}{section}
\newcommand{\spt}{\operatorname{spt}}
\newcommand{\leg}[2]{\left(\frac{#1}{#2}\right)}
\newcommand{\Z}{\mathbb Z}
\newcommand{\F}{\mathbb F}
\newcommand{\CT}{\operatorname{CT}}
\newcommand{\pp}{\operatorname{pp}}
\newcommand{\coeff}[2]{[q^{#1}]#2}
\title[Ramanujan-type congruences for the $\spt$]{Ramanujan-type congruences for the smallest parts function}
\author{Larry Rolen}
\begin{document}
\begin{abstract}
In a famous paper of Andrews, he defined and studied the smallest parts function $\spt(n)$. 
He showed that it satisfies the three Ramanujan-type congruences:
 \[
\spt(5n+4)\equiv0\pmod5,\qquad
\spt(7n+5)\equiv0\pmod7,\qquad
\spt(13n+6)\equiv0\pmod{13}.
\]
In a recent preprint, the author showed that these are the only ``canonical'' congruences.

Here, we complete the classification of Ramanujan-type congruences for $\spt$, showing that if 
$\spt(\ell n+b)\equiv0\pmod\ell$ for all $n\geq0$,
where $\ell\geq5$ is prime, then $(\ell,b)\in\{(5,4),(7,5),(13,6)\}$. That is, Andrews' congruences are the only Ramanujan-type congruences for $\spt(n)$. 
These results have been verified in Lean by AxiomProver.
\end{abstract}
\maketitle

\section{Introduction and statement of results}

The smallest parts function $\spt(n)$ counts the total number of
occurrences of the smallest part in the partitions of $n$. Andrews
introduced this function in \cite{A}, where he proved the congruences
\begin{equation}
 \spt(5n+4)\equiv0\pmod5,\qquad
 \spt(7n+5)\equiv0\pmod7,\qquad
 \spt(13n+6)\equiv0\pmod{13}.
\end{equation}
Many further works on $\spt(n)$  developed the theory of further generalized congruences, crank-type explanations, and arithmetic; see, e.g., 
\cite{AA,ABL,AK,MockModGrids,AGL,Garvan,Ono}.

These are congruences of the simplest type, which we call {\it Ramanujan-type}: congruences modulo a prime $\ell$ of the shape
\begin{equation}\label{eq:simple} 
\spt(\ell n+b)\equiv0\pmod{\ell} \quad \text{ for every } n\geq0.
\end{equation}
All three of Andrews' congruences  satisfy $24b\equiv1\pmod\ell$.
In \cite{R}, the author called such congruences \emph{canonical}. He further classified the canonical congruences by studying
reductions of class polynomials on $X_0(6)$, showing in Theorem~1.1 loc. cit. that Andrews' congruences are the only canonical ones.

Here, we resolve the remaining cases of non-canonical Ramanujan-type congruences.

\begin{theorem}\label{thm:main}
Let $\ell\geq5$ be prime and $0\leq b<\ell$. If
\[
 \spt(\ell n+b)\equiv0\pmod\ell\qquad\text{for every }n\geq0,
\]
then $24b\equiv1\pmod\ell$.
\end{theorem}

Combining Theorem~\ref{thm:main} with the canonical classification gives
the following consequence.

\begin{corollary}\label{cor:classification}
Let $\ell$ be prime and $0\leq b<\ell$. Then \eqref{eq:simple} holds
if and only if
\[
 (\ell,b)\in\{(5,4),(7,5),(13,6)\}.
\]
\end{corollary}
\begin{proof}
These three congruences were established by Andrews.
There are no Ramanujan-type congruences modulo $2$ or $3$, as can be verified using the values:
$\spt(1)=1$, $\spt(2)=3$, $\spt(3)=5$, and $\spt(5)=14$.
When $\ell\geq5$,  Theorem~1.1 of \cite{R} rules out further canonical congruences. Theorem~\ref{thm:main} above rules out 
any non-canonical congruences.
\end{proof}

The first step is to show that any noncanonical congruence propagates to all indices
in one of the two nonzero quadratic-character classes. We then use
constant terms of products of half-integral weight forms to relate
the coefficients of two weakly holomorphic forms. Their lifts to
$X_0(6)$ give a differential congruence. At supersingular points, this
congruence forces a logarithmic derivative to vanish more often than
its degree permits. The same argument treats both character signs.
An explicit character estimate makes this obstruction effective for
large primes, and a finite computation handles the remaining range.

\section*{Notation}
Throughout, $p(n)$ denotes the partition function, $\spt(0):=0$,
$q:=e^{2\pi i\tau}$ for $\tau\in\mathbb H$, and $\Theta:=q\frac{d}{dq}$. We also set
$\chi(m):=\leg{12}{m}$,  and use the
usual normalizations of the modular forms $\eta$, $j$, $E_k$, and $\Delta$.
The Serre derivative, which sends modular forms of weight $k$ to modular forms of weight $k+2$, is $D_k:=\Theta-\frac{k}{12}E_2$.
We denote the multiplier system of $\eta$ on $\operatorname{SL}_2(\Z)$ by $\nu_\eta$. 

\section*{Use of generative AI}
This paper arose out of a conversation between the author and AI. The author takes responsibility for the correctness of the proofs and for the clarity of the exposition. 
The author is also grateful to Connor Olson and Ken Ono for their assistance in autoformalizing the results of this paper.
They used AxiomProver to verify the proofs in Lean. These formal verifications were completely carried out by AxiomProver.
The relevant files are all posted in the following repository:
\url{https://github.com/larryrolen/SPT_noncanonical}

\section{Propagation to a quadratic-character class}\label{sec:propagation}

A congruence for a modular form on one arithmetic progression propagates to
progressions whose indices differ by multiplication by a square.
We need the modularity properties of the generating function of $\spt(n)$. 
This  is closely related to a mock modular form; see \cite{Bringmann}, Section~6 of \cite{Zagier}, and Section~3 of \cite{AA}. 
To show propagation of congruences, we first remove the nonholomorphic part of its modular
completion by a quadratic twist. This yields the following result.

\begin{lemma}\label{lem:prop}
Suppose that \eqref{eq:simple} holds and that $24b\not\equiv1\pmod\ell$.
Let $\sigma:=\leg{1-24b}{\ell}$. Then
\begin{equation}\label{eq:class}
 \spt(N)\equiv0\pmod\ell
 \quad\text{whenever}\quad \leg{1-24N}{\ell}=\sigma.
\end{equation}
\end{lemma}
\begin{proof}
Let $g:=\sum a(n)q^n$ be a weakly holomorphic form of integral or
half-integral weight $w$, level $N_0$, and quadratic character Nebentypus, with
rational $\ell$-integral coefficients. Set
\[
 g_{M,u}:=\sum_n a(Mn+u)q^{n+u/M},\qquad R:=24MN_0.
\]
By \cite[Pages 126--127]{AK15} (see also \cite[Proposition 4]{GJ}), if
$\gamma\in\Gamma_0(24MN_0)$ has upper-left entry $a$ coprime to $6$, then
\[
 (g_{M,u})^R\big|_{Rw}\gamma=(g_{M,a^2u})^R.
\]
Every $a$ coprime to $6MN_0$ can occur as such an entry.

Suppose that $g_{M,u}\equiv0\pmod\ell$. For sufficiently large $s$,
the form $\ell^{-R}\Delta^s(g_{M,u})^R$ is holomorphic and has an
$\ell$-integral expansion. The $q$-expansion principle as in  Corollaire 3.12 of \cite{DR73} 
(see also Corollary~5.3 of \cite{Rad12} and Theorem~5 of \cite{GJ})
shows that integrality is preserved under $\gamma$.
Since $\Delta=q(1+O(q))$ is a unit in the ring of integral Laurent
series, $\ell^{-R}(g_{M,a^2u})^R$ is integral as well. 
As its coefficients
are rational, integrality in the cyclotomic coefficient ring
descends to $\Z_{(\ell)}$. Reducing $(g_{M,a^2u})^R$ modulo $\ell$
in the domain $\F_\ell((q^{1/M}))$ now gives
$g_{M,a^2u}\equiv0\pmod\ell$. In particular, this argument allows
$\ell\mid M$ and $\ell\mid R$.

To apply this observation, define
\[
 \mathcal P:=\eta(24\tau)^{-1},\qquad
 \mathcal S:=\sum_{n\geq0}\spt(n)q^{24n-1},\qquad
 \mathcal M:=\mathcal S+\tfrac1{12}\Theta\mathcal P.
\]
The completion of $\mathcal M$ has weight $3/2$, level $576$, and a
nonholomorphic part supported at negative squares
(see, e.g., \cite[Section 3]{ABL}; the point is that the shadow is a multiple of $\eta(24\tau)$, which can be expressed as a theta function). For a prime $Q\nmid6\ell$, consider the operator
\[
 R_QA:=\tfrac12\bigl(A\otimes\chi_Q^2-\chi_Q(-1)A\otimes\chi_Q\bigr),
 \qquad \chi_Q(u):=\leg{u}{Q},
\]
where $\otimes$ denotes coefficientwise twisting by the character. This operator retains
exactly the indices $u$ for which $\leg{-u}{Q}=-1$, and hence annihilates
the nonholomorphic part. It follows that
$\mathcal M_Q:=R_Q\mathcal M$ and $\mathcal P_Q:=R_Q\mathcal P$ are weakly
holomorphic. Moreover,
\begin{equation}\label{eq:approx}
 E_{\ell-1}\bigl(\mathcal M_Q-\tfrac1{12}D_{-1/2}\mathcal P_Q\bigr)
 +\tfrac1{288}E_{\ell+1}\mathcal P_Q
 \equiv R_Q\mathcal S\pmod\ell.
\end{equation}
The form on the left has weight $\ell+1/2$, level dividing $576Q^4$,
and character $\chi$. The congruence follows from the classical congruences
$E_{\ell-1}\equiv1\pmod\ell$ and $E_{\ell+1}\equiv E_2\pmod\ell$.

Now let $N>0$ satisfy $\leg{1-24N}{\ell}=\sigma$.
We may choose $Q\nmid6\ell$ with $\leg{1-24N}{Q}=-1$.
Indeed, put $D=24N-1$ and $r=24D\ell-1$. Quadratic reciprocity for
Jacobi symbols gives
\[
 \leg{-D}{r}=\leg{r}{D}=\leg{-1}{D}=-1.
\]
Some prime divisor $Q$ of $r$ therefore has the required property,
and $r\equiv-1\pmod{6\ell}$ ensures that $Q\nmid6\ell$.
Choose $a$ with
\[
 24N-1\equiv a^2(24b-1)\pmod\ell,\qquad a\equiv1\pmod{24Q}.
\]
The Chinese Remainder Theorem supplies such an $a$, and it is coprime
to $6\ell Q$. Since $a^2\equiv1\pmod{24}$, the first congruence also
holds modulo $24\ell$. Applying the propagation argument to
\eqref{eq:approx}, with $M=24\ell$, proves that its coefficient at
$q^{24N-1}$ vanishes modulo $\ell$. This coefficient is $\spt(N)$
modulo $\ell$. Finally, $N=0$ causes no difficulty since $\spt(0)=0$.
\end{proof}

\section{Two constant-term identities}\label{sec:constantterms}

We next translate \eqref{eq:class} into a relation between coefficients
of weakly holomorphic forms. It is convenient to use the unscaled
functions
\[
 P:=\eta^{-1},\qquad S:=\sum_{N\geq0}\spt(N)q^{N-1/24},\qquad
 f:=-12S-24\Theta P=\sum_r a(r)q^{r/24}.
\]
Thus $a(-1)=1$ and
$a(24N-1)=-12\spt(N)-(24N-1)p(N)$. The harmonic Maass form completion of $f$ has
weight $3/2$, multiplier $\overline{\nu_\eta}$, and shadow proportional
to $\eta$.

For $d>0$ with $d\equiv23\pmod{24}$, consider the integral weakly
holomorphic forms
\begin{align}
 H_d&:=q^{-d/24}-a(d)q^{1/24}+O(q^{25/24})
       \in M^!_{1/2}(\nu_\eta),\label{eq:H}\\
 K_d&:=q^{-d/24}+\sum_{r>0}k_d(r)q^{r/24}
       \in M^!_{5/2}(\nu_\eta),\label{eq:K}
\end{align}
and write $H_d=q^{-d/24}+\sum_{r>0}h_d(r)q^{r/24}$.
These forms can be  obtained by row reduction on principal parts in $\eta\Z[j]$
and $-\eta\Theta j\,\Z[j]$, respectively. Using the notation $\CT$ to denote the
constant term, \eqref{eq:H} shows that
$\CT(H_df)=0$. This gives the first of the two identities below;
the second comes from the product of $K_d$ with $P$.

\begin{lemma}\label{lem:coeff}
Assume \eqref{eq:class}, and let $d>0$ satisfy
$d\equiv23\pmod{24}$ and $\leg{-d}{\ell}=\sigma$. Then
\begin{equation}\label{eq:squares}
 m^2h_d(m^2)+d\,k_d(m^2)\equiv0\pmod\ell
 \qquad\text{for all }m\geq1\text{ with }(m,6\ell)=1.
\end{equation}
\end{lemma}
\begin{proof}
For a prime $p\nmid6$, let $T_p^{(w)}$ denote the Hecke operator of
index $p^2$. In weights $w=-1/2$ and $3/2$, it acts on Fourier expansions by:
\[
 (A|T_p^{(w)})(r)=A(p^2r)+\chi(p)\leg{-r}{p}p^{w-3/2}A(r)
                       +p^{2w-2}A(r/p^2).
\]
Here a coefficient with a nonintegral index is interpreted as zero.
Factoring $m=p_1\cdots p_s$, with repeated prime factors allowed, set
$U_m^{(w)}:=T_{p_1}^{(w)}\cdots T_{p_s}^{(w)}$ and $U_1^{(w)}:=1$.
We use these compositions throughout. Comparing coefficients gives
\[
 (\Theta P)|U_m^{(3/2)}=m^2\Theta(P|U_m^{(-1/2)}).
\]
Let $\pp$ denote the principal part of a $q$-series. Since a negative
square index $-t^2$ can only change to $-(pt)^2$, $-t^2$, or
$-(t/p)^2$, there are integers $c_t$ such that
\begin{align*}
 \pp(f|U_m^{(3/2)})&=\sum_{t\mid m}c_tq^{-t^2/24},\qquad c_m=m,\\
 \pp(m^2P|U_m^{(-1/2)})&=\sum_{t\mid m}\frac{c_t}{t^2}q^{-t^2/24}.
\end{align*}
The second identity also follows directly from the derivative relation
and $\pp(f)=\pp(-24\Theta P)$.

Iterating \cite[Lemma 2.4]{R}, we find that
\[
 f|U_m^{(3/2)}-\lambda_m f\in M^!_{3/2}(\overline{\nu_\eta}), \text{ where} 
 \qquad \lambda_m:=\prod_{i=1}^s\chi(p_i)(1+p_i).
\]
A weakly holomorphic modular form of weight $2$ and level one has
constant term zero.
Apply this first to $H_d$
times the preceding difference, and then to $K_d(P|U_m^{(-1/2)})$.
Using $\CT(H_df)=0$, we obtain
\begin{align}
 \coeff{d/24}{(f|U_m^{(3/2)})}
       &=-\sum_{t\mid m}c_t h_d(t^2),\label{eq:ct1}\\
 m^2\coeff{d/24}{(P|U_m^{(-1/2)})}
       &=-\sum_{t\mid m}\frac{c_t}{t^2}k_d(t^2).\label{eq:ct2}
\end{align}

On the character class in \eqref{eq:class}, we have
$f\equiv-24\Theta P\pmod\ell$, including at negative indices.
Each prime Hecke operator occurring here preserves that class.
Consequently, the left-hand side of \eqref{eq:ct1} is $-d$ times the
left-hand side of \eqref{eq:ct2} modulo $\ell$. We conclude that
\[
 \sum_{t\mid m}\frac{c_t}{t^2}
       \bigl(t^2h_d(t^2)+d\,k_d(t^2)\bigr)\equiv0\pmod\ell.
\]
The result now follows by (strong) induction on $m$. All proper-divisor terms
vanish by induction, and the remaining coefficient is
$c_m/m^2=1/m$, which is an $\ell$-adic unit.
\end{proof}

\section{The obstruction on \texorpdfstring{$X_0(6)$}{X0(6)}}\label{sec:geometry}

The coefficient relation in Lemma~\ref{lem:coeff} becomes particularly
useful after lifting $H_d$ and $K_d$ to integral weight. We recall the
normalizations of these lifts before giving the obstruction.

Suppose that $d$ is squarefree and $d\equiv23\pmod{24}$, and put
$h:=h(-d)$. For $r\in\{1,5,7,11\}$, consider the $\Gamma_0(6)$-classes
of primitive positive forms $[a,b,c]$ of discriminant $-d$ with
$6\mid a$ and $b\equiv r\pmod{12}$. The corresponding CM points are
$\tau=(-b+i\sqrt d)/(2a)$. The function
\[
 t(\tau)=\frac{\eta(\tau)^5\eta(3\tau)}{\eta(2\tau)\eta(6\tau)^5}
\]
is a coordinate on $X_0(6)_{\Z[1/6]}\cong\mathbb P^1$.
Let $A_d(X)$ be the product of $X-t(\tau)$ over the classes with
$r=1,11$, and let $B_d(X)$ be the corresponding product for $r=5,7$.
Both are monic polynomials of degree $2h$ in $\Z[1/6][X]$.
The normalized Borcherds product of \cite[Section 4]{R} is
\begin{equation}\label{eq:product}
 \Psi_d:=\frac{A_d(t)}{B_d(t)}
       =\prod_{m\geq1}(1-q^m)^{\chi(m)h_d(m^2)}.
\end{equation}
Taking its logarithmic derivative gives a meromorphic form of weight $2$:
\begin{equation}\label{eq:L}
 L_d:=\Theta\log\Psi_d
     =-\sum_{m\geq1}\chi(m)m h_d(m^2)\frac{q^m}{1-q^m}.
\end{equation}
We will use Borcherds' work \cite[Theorem 14.3]{B} to build a weight $4$
 meromorphic Shimura lift of $K_d$:
\begin{equation}\label{eq:J}
 J_d=\sum_{m\geq1}\chi(m)k_d(m^2)\frac{q^m}{(1-q^m)^2}
\end{equation}
on $\Gamma_0(6)$. This is holomorphic at the cusps and has poles of
order at most two, all at the CM points above.

To specify the Shimura lift input, realize $\Lambda=U\oplus\langle12\rangle$
as $\Z^3$ with quadratic form $Q(a,b,c)=6b^2-ac$. Its dual has
$b\in\frac1{12}\Z$, and its discriminant form is
$(\Z/12\Z,r\mapsto r^2/24)$. The vector
\[
 v=\mathfrak e_1-\mathfrak e_5-\mathfrak e_7+\mathfrak e_{11}
\]
transforms by $\nu_\eta$ under the Weil representation. We apply \cite[Theorem 14.3]{B} to $K_dv$
with $m^+=2$, noting that
$\Gamma_0(6)$ lies in the discriminant kernel. At infinity, the positive
dual vectors have norm $m^2/24$ and pair with the coordinate as
$m\tau$. The Fourier expansion formula in \cite[Theorem 5.1]{LZ} is therefore
$\sum_{u,m\geq1}u\chi(m)k_d(m^2)q^{um}$, with zero constant term,
which is \eqref{eq:J}. A negative vector $(a,r/12,c)$ of norm $-d/24$
gives the quadratic $6c\tau^2+r\tau+a$ of discriminant $-d$.
This identifies the poles. Finally, the Atkin--Lehner involutions
permute the four cusps and act on $v$ by signs, so the vanishing at
infinity gives holomorphy at all cusps.

We also need the following separation of the two sets of CM points.
The point here is that the same bound applies to both splitting types.

\begin{lemma}\label{lem:separation}
If $\ell\nmid6d$ and $\ell>(d+1)^2/24$, then
\begin{equation}\label{eq:sep}
 (\overline A_d,\overline B_d)=1\quad\text{in }\F_\ell[X].
\end{equation}
\end{lemma}
\begin{proof}
A common root gives a reduced pair $(E,C_6)$ (here, $C_6$ is a cyclic group of order $6$) with two CM actions
$\alpha,\beta$ of $(1+\sqrt{-d})/2$, each of trace $1$ and norm
$(d+1)/4$. One acts on $C_6$ by $0$ or $1$, and the other by $3$ or $4$.
If they commuted, then $\beta=\alpha$ or $1-\alpha$, which gives a contradiction.
 In particular, ordinary reduction is impossible.

In the supersingular case, $\operatorname{End}(E,C_6)$ is an Eichler
order of discriminant $6\ell$. The independent elements
$1,\alpha,\beta,\alpha\beta$ span a sublattice whose covolume for the
trace norm is at least $6\ell$. On the other hand, Hadamard's inequality
bounds this covolume by
\[
 4\operatorname N(\alpha)\operatorname N(\beta)=\frac{(d+1)^2}{4}<6\ell.
\]
This contradiction proves the claim (see also \cite[Lemma 5.2]{R}).
\end{proof}

We can now state the obstruction in the form needed for the proof.

\begin{proposition}\label{prop:obstruction}
Suppose that there is a squarefree $d>0$ with $d\equiv23\pmod{24}$,
$\leg{-d}{\ell}=\sigma$, and
\begin{equation}\label{eq:bounds}
 \ell>\frac{(d+1)^2}{24},\qquad \ell+1>8h(-d).
\end{equation}
Then \eqref{eq:class} cannot hold.
\end{proposition}
\begin{proof}
Suppose that \eqref{eq:class} holds. Differentiating \eqref{eq:L} once
and applying Lemma~\ref{lem:coeff}, we see that the terms with
$\ell\nmid m$ agree with those of $dJ_d$ modulo $\ell$.
A second derivative annihilates the remaining terms, since their
exponents are divisible by $\ell$. Hence,
\begin{equation}\label{eq:differential}
 \Theta^2L_d\equiv d\,\Theta J_d\pmod\ell.
\end{equation}

Let $H:=\overline E_{\ell-1}$ be the Hasse invariant and put
$E=\overline E_{\ell+1}$. The operator
\[
 \vartheta_kF:=H D_kF+\frac{k}{12}EF
\]
has weight $k+\ell+1$ and $q$-expansion $\Theta F$ modulo $\ell$.
Thus, after making the weights equal, \eqref{eq:differential} gives
\begin{equation}\label{eq:hasse}
 \vartheta_{\ell+3}\vartheta_2 L_d=dH\vartheta_4J_d.
\end{equation}
Here, reduction away from the CM support follows from the
$q$-expansion principle at good reduction. More explicitly,
$(A_d(t)B_d(t))^2\Delta^{8h}J_d$ has no poles and has integral
Fourier coefficients. The same pole-clearing argument applies to
the derivatives. For $L_d$, one can also use the rational expression
below.

The zeros of $H$ are simple and are precisely the supersingular points
\cite[Sections 2 and 4.3]{Katz}. Since $E=12D_{\ell-1}H$, the form
$E$ is nonzero at each such point. Evaluating \eqref{eq:hasse} at a
supersingular point outside the CM support gives
$E^2L_d/24=0$. Thus, $L_d$ vanishes at every one of these points.
The curve $X_0(6)$ has genus zero, four cusps, and no elliptic points,
so its weight-$2$ bundle has degree $2$. It follows that $H$ has degree
$\ell-1$, and hence that there are $\ell-1$ supersingular points over
$\overline{\F}_\ell$. At most $4h$ lie in the CM support.

On the other hand, logarithmic differentiation of \eqref{eq:product}
gives
\[
 L_d=(\Theta t)\frac{W_d(t)}{A_d(t)B_d(t)},\qquad
 W_d=A_d'B_d-A_dB_d',\qquad \deg W_d\leq4h-2.
\]
Lemma~\ref{lem:separation} and the inequality $2h<\ell$ imply that
$\overline W_d\ne0$. Indeed, otherwise $\overline A_d$ would divide
$\overline A_d'$, forcing $\overline A_d'=0$. This is impossible
for a monic polynomial of degree $2h<\ell$.
Since $t$ is a coordinate, the Kodaira--Spencer isomorphism shows that
$\Theta t$ is nonvanishing on the interior. Therefore $\overline W_d$
has at least
\[
 \ell-1-4h>4h-2
\]
distinct roots. This contradicts its degree bound.
\end{proof}

\section{Proof of Theorem~\ref{thm:main}}\label{sec:proof}

It remains to produce the discriminants required by
Proposition~\ref{prop:obstruction} for large primes, and to exclude
\eqref{eq:class} directly for the smaller primes.

\subsection{Large primes}
To handle large primes, we require a bound of Trevi\~no.
Using the correction in
\cite[Theorem~1]{T}, his Theorem~1 implies that if a nonprincipal character modulo
an odd prime $\ell$ is constant on $H$ consecutive integers, then
\begin{equation}\label{eq:character}
H\le3.38\cdot \ell^{1/4}\log \ell.
\end{equation}
Fix $\sigma\in\{1,-1\}$, and let $k\geq1$ be least with
$\leg{1-24k}{\ell}\ne-\sigma$. After translating by
$24^{-1}\pmod\ell$, this estimate gives
$k\leq3.38\ell^{1/4}\log\ell+1$.
For $\ell\geq10^{11}$, set $d:=24k-1$. Then we have
\[
 d+1\leq81.12\ell^{1/4}\log\ell+24<\sqrt{24\ell}.
\]
For the last inequality, the ratio of the middle expression to
$\sqrt{24\ell}$ is less than $0.785$ at $10^{11}$ and decreases
thereafter. In particular, $d<\ell$, so $\leg{-d}{\ell}=\sigma$.

The minimality of $k$ implies that $d$ is squarefree. Indeed, if
$u^2\mid d$ with $u>1$, then $d/u^2$ still belongs to the residue class
$23$ modulo $24$ and has the same character value, giving a smaller
choice of $k$. Counting reduced forms gives
\[
 h(-d)\leq\frac d3+2\sqrt{\frac d3}<d.
\]
Consequently,
$8h(-d)<8\sqrt{24\ell}<\ell$. Both inequalities in
\eqref{eq:bounds} hold, so Proposition~\ref{prop:obstruction} excludes
\eqref{eq:class} for every $\ell\geq10^{11}$.

\subsection{The finite range}
For $5\leq\ell<10^{11}$, we computed a nonzero coefficient in each
character class. Let
\[
 \mathcal Q:=\{q\leq2087:q\text{ prime},\ q\equiv23\pmod{24}\},
 \qquad s_q:=\spt((q+1)/24).
\]
This gives $40$ values $s_q$, with indices at most $87$.
Let $\mathcal E$ be the set of primes at least $5$ dividing
$\prod_{q\in\mathcal Q}q s_q$. Exact factorization gives $84$ such
primes. For each $\ell\in\mathcal E$, the accompanying certificate
gives $N_+,N_-\leq8$ such that
\[
 \ell\nmid\spt(N_\pm),\qquad \leg{1-24N_\pm}{\ell}=\pm1.
\]
The $\spt$-values were computed using the generating function 
\[\sum_{n\geq0}\spt(n)q^n
 =\sum_{m\geq1}\frac{q^m}{(1-q^m)^2\prod_{j>m}(1-q^j)}
\]
from the proof of Theorem~4 in \cite{A}.

If $\ell\notin\mathcal E$, every $s_q$ is nonzero modulo $\ell$, and
quadratic reciprocity gives $\leg{-q}{\ell}=\leg{\ell}{q}$.
It therefore suffices to verify the finite claim
\begin{equation}\label{eq:finite}
 \begin{gathered}
 5\leq x\leq10^{11},\quad (x,6\prod_{q\in\mathcal Q}q)=1,\quad
 \leg{x}{q}\text{ independent of }q\in\mathcal Q\\
 \Longrightarrow\quad x\text{ is a square}.
 \end{gathered}
\end{equation}
A prime cannot be a square, so \eqref{eq:finite} supplies both signs
among the forty symbols.
We verified \eqref{eq:finite} using SAGE \cite{sage}; this took 36.8 seconds on the author's Macbook M5 pro.

\begin{proof}[Proof of Theorem~\ref{thm:main}]
If a noncanonical simple congruence existed, Lemma~\ref{lem:prop}
would give \eqref{eq:class} for one of the two signs $\sigma$.
The large-prime argument excludes this when $\ell\geq10^{11}$,
and the finite computation excludes it for $5\leq\ell<10^{11}$.
\end{proof}

\end{document}